\documentclass[review,12pt,authoryear]{elsarticle}

\usepackage{multirow}
\usepackage{amsmath,xfrac,cases,mathrsfs} 
\usepackage{amssymb}
\usepackage{enumerate}
\usepackage{lineno,hyperref}
\usepackage{braket}

\newtheorem{thm}{Theorem}
\newtheorem{lem}[thm]{Lemma}
\newdefinition{rmk}{Remark}
\newdefinition{prop}[thm]{Proposition}
\newproof{pf}{Proof}
\newproof{pot}{Proof of Theorem \ref{thm2}}

\journal{!}

\begin{document}
	
\begin{frontmatter}
		
\title{Nonparametric Estimation for Stochastic Differential Equations Driven by Fractional Brownian Motion}

\author[rvt]{Han Yuecai}
\ead{hanyc@jlu.edu.cn}
		
\author[rvt]{Zhang Dingwen\corref{cor1}}
\ead{zhangdw20@mails.jlu.edu.cn}

\cortext[cor1]{Corresponding author}

\address[rvt]{School of Mathematics, Jilin University, Changchun, 130000, China}
		
\begin{abstract}
We study the nonparametric Nadaraya-Watson estimator of the drift function for ergodic stochastic processes driven by fractional Brownian motion of Hurst parameter $H>1/2$. The estimator is based on the discretely observed stochastic processes. By using the ergodic properties and stochastic integral, we obtain the consistency of the proposed estimator. 
\end{abstract}
		
		
\end{frontmatter}
	
	
\section{Introduction}
We consider the following stochastic differential equation (SDE in short) driven by fractional Brownian motion (fBm in short)
\begin{equation}{\label{eq1}}
\left\{
\begin{aligned}
\mathrm{d}X_{t}&=b(X_{t})\mathrm{d}t+\sigma\mathrm{d}B^{H}_{t},\\
X_{0}&=x,
\end{aligned}
\right.
\end{equation}
where $x\in\mathbb{R}$ is the initial value of the process $\{X_{t}, t\geq 0\}$, $b:\mathbb{R}\rightarrow\mathbb{R}$ is an unknown continuous function, $\sigma\in \mathbb{R}$ is a constant and $\{B^{H}_{t}, t\geq 0\}$ is a fBm of Hurst index $H\in(\frac{1}{2},1)$ defined on a complete probability space $(\Omega, \mathcal{F}, \{\mathcal{F}_{t}\}_{t\geq0}, \mathbb{P})$, where $\mathcal{F}_{t}$ is the $\sigma$-algebra generated by $\{B_{s}^{H}, s\leq t\}$. We consider the stochastic process $\{X_{t}\}$ observed at regularly
sapced time points $\{t_{k}=k\alpha_{n}, k=0, 1, \cdots, n\}$, where $\alpha_{n}$ is the time frequency and $n$ is the sample size. Based on the high-frequency observations, we are aiming to give an estimation for the drift function $b(\cdot)$.

With the development of technology, the parameter estimation problem in SDEs has gained much attention in recent years due to their increased applications in the broad fields. It is necessary that the parameters which characterize the
SDE system should be estimated via the data in many real world applications. 

Along the last decades, there is a huge literture devoted to the problems of parameter estimation. When the drit term is linear which means that $b(X_{t})=-\theta X_{t}$, then $X$ is a fractional Ornstein–Uhlenbeck processes and the estimation of $\theta$ has been widely studied. \cite{Le Breton(2002)} propose a maximum likehood estimator via the approximated likelihood ratio based on continuous or discrete observations. On the other hand, \cite{Hu(2010)} propose the least squares estimator, they consider the integral with respect to $B^{H}$ for $H\in[\frac{1}{2},\frac{3}{4})$. They extend the former work and established a central limit theorem for the least squares estimator for $H\in(0,\frac{3}{4}]$ and a noncentral limit theorem for $H\in(\frac{3}{4},1)$ \citep{Hu(2019)}. Moreover, a ergodic estimator for $\theta$ is also proposed \citep{Hu(2010),Hu(2019)}. When the drift term is $\theta b(X_{t})$, the maximum likehood estimator of $\theta$ is studied both with continuous and discrete observations \citep{Tudor(2007)}. \cite{Hu(2018)} propose the least squares estimator for $\theta$ with $H\in (\frac{1}{4},1)$, they also derive a maximum inequality for It$\hat{o}$-Skorohod integrals. When the drift term is $b(\theta, X_{t})$, \cite{Neuenkirch(2014)} propose a least squares estimator for $\theta$ based on discrete observations of the process $\{X_{t}, t\geq 0\}$.

In reality, the drift function is seldom known, hence it is significant and needed urgently to give a reasonable estimator. On the nonparametric estimation of drift in Eq. (\ref{eq1}), there are only a few references. The most popular used tools are kernel function $K_{h}(\cdot)$ where $h$ is the bandwidth.
In \cite{Saussereau(2014)}, a Nadaraya–Watson estimator is defined at time $t_{n}$ as
\begin{equation}\label{eq2}
\hat{b}_{t_{n} . h}(x)=\frac{\sum_{k=0}^{n-1}\left(t_{n}-t_{k}\right)^{1-2 H} K\left(\left(X_{t_{k}}-x\right) / h\right)\left(X_{t_{k+1}}-X_{t_{k}}\right)}{\sum_{k=0}^{n-1}\left(t_{n}-t_{k}\right)^{1-2 H} K\left(\left(X_{t_{k}}-x\right) / h\right)\left(t_{k+1}-t_{k}\right)}.
\end{equation}
They obtain the consistency of the estimator under a one sided dissipative Lipschitz condition which insures the ergodic property for the solution of Eq. (\ref{eq1}). 
In \cite{Fabienne(2019)}, they propose a Nadaraya–Watson estimator as 
\begin{equation*}
\widehat{b}_{T, h}(x):=\frac{\int_{0}^{T} K\left(\frac{X(s)-x}{h}\right) \mathrm{d} X(s)}{\int_{0}^{T} K\left(\frac{X(s)-x}{h}\right) \mathrm{d} s},
\end{equation*}
where the stochastic integral with respect to $X$ is taken in the sense of It$\hat{o}$-Skorokhod, which is similar to the Eq. (\ref{eq2}).
In \cite{Prakasa(2011)}, they propose an estimator of trend function $b_{t}=b(x^{0}_{t})$ as following
\begin{equation}\label{eq3}
\hat{b}_{t}=\frac{1}{h} \int_{0}^{T} K\left(\frac{\tau-t}{h}\right) \mathrm{d} X_{\tau},
\end{equation}
where $x^{0}_{t}$ is the solution of Eq. (\ref{eq1}) when $\sigma=0$. They obtain the asymptotic behaviour of the estimator Eq. (\ref{eq3}) when $\sigma\rightarrow0$.
One can refer to \cite{Kutoyants(2004)} for nonparameter estimation in diffusion processes.

As we all know, for a nonparameter $b(\cdot)$, one can usually give the estimator \citep{N(1964),W(1964)}
\begin{equation}\label{eq4}
\hat{b}(x)=\int_{-\infty}^{\infty}K_{h}(y-x)b(y)\mathrm{d}y.
\end{equation}
Motivated by the aforementioned works, in this paper, we propose a Nadaraya-Waston estimator of the drift function $b(\cdot)$ based on the discrete obervations as 
\begin{equation}\label{eq5}
\hat{b}_{n,h}(x)=\sum_{k=0}^{n-1}W_{n,k}(X_{t_{k}}, x)\diamond(\frac{X_{t_{k+1}}-X_{t_{k}}}{\alpha_{n}}),
\end{equation}
where `$\diamond$' is the Wick product which has been associated with integrals of It$\hat{o}$-Skorokhod type \citep{Hu(1996),Holden(1996)} and the weight function $W_{n,k}(X_{t_{k}},x)$ is given by
\begin{equation*}
W_{n,k}(X_{t_{k}},x)=\frac{K_{h}(X_{t_{k}}-x)}{\sum_{k=0}^{n-1}K_{h}(X_{t_{k}}-x)},
\end{equation*}
where $K_{h}(\cdot)=1/hK(\cdot/h)$ for $k=0, 1, \cdots, n$, $K(\cdot)$ is a kernel function and $h$ is the bandwidth. The estimator defined in Eq. (\ref{eq5}) is an analogue of Eq. (\ref{eq4}).

This paper is organized as follows. In Section \ref{sec2}, we describe some priliminary results on stochastic integral. Based on some assumptions, we obtain the ergodic properties of $X_{t}$ and a bound for the Malliavin derivative of $X_{t}$. Our main results are proposed in Section \ref{sec3}. Based on some prepared lemmas, we obtain the consistency of the estimator.

\section{Preliminaries}\label{sec2}
In this section, we give some notations and state our assumptions.

The fBm $\left(B_{t}^{H}, t \in \mathbb{R}\right)$ with Hurst parameter $H \in(0,1)$ is a zero mean Gaussian process with covariance
\begin{equation*}
\mathbb{E}\left(B_{t}^{H} B_{s}^{H}\right)=R_{H}(s, t)=\frac{1}{2}\left(|t|^{2 H}+|s|^{2 H}-|t-s|^{2 H}\right).
\end{equation*}
On any finite interval, all paths of fBm are $\alpha$-H$\ddot{o}$lder  continuous with $\alpha<H$. Denote by $\eta_{T}$ the coefficient of fBm on the interval $[0,T]$, it means that
\begin{equation}\label{eq6}
\sup_{t\neq s\in [0,T]}|B^{H}_{t}-B^{H}_{s}|\leq\eta_{T}|t-s|^{\alpha}.
\end{equation}
Futhermore, we have $\mathbb{E}|\eta_{T}|^{p}=T^{p(H-\alpha)}\mathbb{E}|\eta_{1}|^{p}$ by the self-similarity property of fBm for any $p>1$. Let $f ,g :\mathbb{R}_{+}\rightarrow \mathbb{R}$ be Borel measurable functions. Denote by $\mathcal{H}$ the Hilbert space equipped with the inner product
\begin{equation*}
\Braket{f,g}_{\mathcal{H}}=H(2H-1)\int_{0}^{\infty}\int_{0}^{\infty}f(t)g(s)|t-s|^{2H-2}\mathrm{d}t\mathrm{d}s.
\end{equation*}
Then $\mathcal{H}$ is a Banach space with the norm $\|\cdot\|_{\mathcal{H}}$. Futhermore, for any $f\in L^{\frac{1}{H}}([0,\infty); \mathbb{R})$, we have 
\begin{equation}\label{eq7}
\|f\|_{\mathcal{H}}\leq c_{H}\|f\|_{L^{\frac{1}{H}}},
\end{equation}
where $c_{H}$ is a constant with respect to $H$ defined in \cite{Nualart(2006)}.

Let $B_{H}(\phi)=\int_{0}^{\infty}\phi_{t}\mathrm{d}B_{t}^{H}$. For a smooth and cylindrical random variable $F=f(B^{H}(\phi_{1}),B^{H}(\phi_{2}),\cdots,B^{H}(\phi_{d}))$, we define its Malliavin derivative as the $\mathcal{H}$-valued random variable as
\begin{equation*}
\mathrm{D}F=\sum_{i=1}^{d}\frac{\partial f}{\partial x_{i}}\left(B^{H}(\phi_{1}),B^{H}(\phi_{2}),\cdots,B^{H}(\phi_{d})\right)\phi_{i}.
\end{equation*}

Let $f\in L^{p}(\Omega, \mathcal{F}, \mathbb{P})$ for each $p\geq 1$, the It$\hat{o}$ type stochastic integral is the limit of the Riemann sums defined in terms of the Wick product \citep{Duncan(2000)}
\begin{equation}\label{eq8}
\int_{0}^{T}f_{s}\mathrm{d}B^{H}_{s}=\lim_{\alpha_{n}\rightarrow 0}\sum_{k=0}^{n-1}f_{t_{k}}\diamond(B^{H}_{t_{k+1}}-B^{H}_{t_{k}}).
\end{equation}

We make use of the notation $\delta(u)=\int_{0}^{\infty}u_{t}\mathrm{d}B^{H}_{t}$ and
and call $\delta(u)$ the divergence integral of $u$ with respect to the fBm $B_{\cdot}^{H}$. The relationship between divergence integral and pathwise integral is stated as the following formula
\begin{equation*}
\int_{0}^{T}u_{t}\mathrm{d}B_{t}^{H}=\int_{0}^{T}u_{t}\delta B_{t}^{H}+H(2H-1)\int_{0}^{T}\int_{0}^{t}\mathrm{D}_{s}u_{t}|t-s|^{2H-2}\mathrm{d}s\mathrm{d}t,
\end{equation*}
where `$\delta$' means that the integral is taken in the sense of Stratonovich.
Let $\|f\|_{\infty}=\sup_{\cdot\in \mathbb{R}} f(\cdot)$.

We will make use of the following assumptions:
\begin{enumerate}[i.]
	\item\label{asp1} The dift function $b(\cdot)$ satisfies a global Lipschitz conditions which is there exists a positive constant $L>0$, such that 
	\begin{equation*}
		|b(x)-b(y)|\leq L|x-y|, \quad x, y \in \mathbb{R}.
	\end{equation*}
	\item\label{asp2} The dift term $b(\cdot)$ is continuously differentiable with a polynomial growth condition on its derivative $b^{\prime}(\cdot)$ and $b(\cdot)$ itself, there exists a constant $C>0$ and $m\in \mathbb{N}$ such that
	\begin{equation*}
		|b(x)|+|b^{\prime}(x)|\leq C(|1+|x|^{m}), x\in \mathbb{R}.
	\end{equation*}
	\item\label{asp3} There exists a constant $M>0$ such that
	\begin{equation*}
		\left(b(x)-b(y)\right)\times(x-y)\leq -M|x-y|^{2}, x, y\in \mathbb{R}. 
	\end{equation*} 
	\item\label{asp4} The observation time frequency $\alpha_{n}=O(n^{-1+\frac{1}{\gamma}})$, $t_{n}=n\alpha_{n}\rightarrow \infty$ as $n\rightarrow \infty$ where $\gamma>\max\{1+m^{2}H,2\}$.
	\item\label{asp5} The kernel function $K(\cdot)$ is continuously differentiable, non-negative with support $[-1, 1]$. The bandwidth $h$ satifies $h\rightarrow 0$ as $n\rightarrow \infty$.
	\item\label{asp6} $\|b(\cdot)\|_{\infty}+\|b^{\prime}(\cdot)\|_{\infty}<\infty$ and $\|K\|_{\infty}<\infty$.
\end{enumerate}
We define the shift operators $\theta_{t}: \Omega\rightarrow\Omega$ as 
\begin{equation*}
\theta_{t}w(\cdot)=w(\cdot+t)-w(t), t\in\mathbb{R}, w\in\Omega.
\end{equation*}
We summarizes the ergodic property for the process $\{X_{t}, t\geq0\}$ in the following theorem, and for a proof one can see in \cite{Saussereau(2014)}. One can see in \cite{Hairer(2005),Hairer(2007),Haier(2011)} for more about the ergodic properties of diffusion processes.
\begin{thm}\label{th1}
Under the conditions \ref{asp2}, \ref{asp3} and \ref{asp4}. Let $\varphi$ be a continuously differentiable function such that
\begin{equation*}
|\varphi(x)|+|\varphi^{\prime}(x)|\leq C_{\varphi}|1+|x|^{p}|, \quad y \in \mathbb{R}, 
\end{equation*}
for some $C_{\varphi}>0$ and $p\in \mathbb{N}$. There extists a random variable $\bar{X}$ with finite moments of any order such that 
\begin{enumerate}
	\item We have
	\begin{equation*}
		\lim _{T \rightarrow \infty} \frac{1}{T} \int_{0}^{T} \varphi\left(X_{s}\right) \mathrm{d} s=\mathbb{E}(\varphi(\bar{X})),
	\end{equation*}
    almost surely.
	\item If $\gamma>1+\left(m^{2}+p\right) H$ and $\gamma>p+1$ then
	\begin{equation}\label{eq9}
	\lim _{n \rightarrow \infty} \frac{1}{t_{n}} \int_{0}^{t_{n}}\left\{\sum_{k=0}^{n-1} \varphi\left(X_{t_{k}}\right) \mathbf{1}_{\left[t_{k}, t_{k+1}\right)}(s)\right\} \mathrm{d} s=\mathbb{E}(\varphi(\bar{X})),
	\end{equation}	
    almost surely.
\end{enumerate}
\end{thm}
Next, we give a upper bound for the Malliavin derivative of $X_{t}$. More estimates for the p-th moment of $X_{t}$ one can refer to Proposition 4.1 in \cite{Hu(2019)}.
\begin{prop}\label{prop1}
Under the contions \ref{asp1}, \ref{asp2} and \ref{asp3}, we have that the Malliavin derivative of $X_{t}$ defined in Eq. (\ref{eq1}) satifies that for all $0\leq s\leq t$,
\begin{equation*}
|\mathrm{D}_{s}X_{t}|\leq |\sigma|e^{-M(t-s)}.
\end{equation*}
\end{prop}
\noindent\textbf{Proof of Proposition \ref{prop1}}.
Note that
\begin{equation*}
\mathrm{D}_{s}X_{t}=\int_{s}^{t}b^{\prime}(X_{u})D_{s}X_{u}\mathrm{d}u+\sigma.
\end{equation*}
Denote $Z_{t}=\mathrm{D}_{s}X_{t}$ for $t>s$. We can write the above equation as
the following ordinary differential equation for $t \geq s$:
\begin{equation*}
\mathrm{d}Z_{t}=b^{\prime}(X_{t})Z_{t}\mathrm{d}t, \quad\quad Z_{s}=\sigma.
\end{equation*}
The solution of the above equation is 
\begin{equation*}
	\mathrm{D}_{s}X_{t}=\sigma e^{\int_{s}^{t}b^{\prime}(X_{r})\mathrm{d}r}.
\end{equation*}
Moreover, differentiating $|Z_{t}|^{2}$ with respect to $t$ and by \ref{asp3}, we have 
\begin{equation*}
\frac{\mathrm{d}|Z_{t}|^{2}}{\mathrm{d}t}=2 b^{\prime}(X_{t})Z^{2}_{t}\leq -2M|Z_{t}|^{2}.
\end{equation*}
Then we obtain the desired results directly by Gronwall's inequality.
\hfill$\square$
\section{Asymptotic behavior of the estimator}\label{sec3}
In this section, we consider asymptotic behavior of the Nadaraya-Watson estimator defined in Eq. (\ref{eq5}) of the drift function. 

First note that
\begin{equation*}
X_{t_{k+1}}-X_{t_{k}}=\int_{t_{k}}^{t_{k+1}}\left(b(X_{s})-b(X_{t_{k}})\right)\mathrm{d}s+b(X_{t_{k}})\alpha_{n}+\sigma(B_{t_{k+1}}^{H}-B_{t_{k}}^{H}),	
\end{equation*}
then the N-W estimator defined in Eq. (\ref{eq5}) can be partitioned into three components as follows:
\begin{equation}\label{eq10}
\begin{aligned}
\hat{b}_{n,h}(x)
=&\frac{\frac{1}{n\alpha_{n}}\sum_{k=0}^{n-1}K_{h}(X_{t_{k}}-x)\diamond(X_{t_{k+1}}-X_{t_{k}})}{\frac{1}{n}\sum_{k=0}^{n-1}K_{h}(X_{t_{k}}-x)}\\
=&\frac{\frac{1}{n\alpha_{n}}\sum_{k=0}^{n-1}K_{h}(X_{t_{k}}-x)\left(\int_{t_{k}}^{t_{k+1}}\left(b(X_{s})-b(X_{t_{k}})\mathrm{d}s\right)\right)}{\frac{1}{n}\sum_{k=0}^{n-1}K_{h}(X_{t_{k}}-x)}\\
&+\frac{\frac{1}{n}\sum_{k=0}^{n-1}K_{h}(X_{t_{k}}-x)b(X_{t_{k}})}{\frac{1}{n}\sum_{k=0}^{n-1}K_{h}(X_{t_{k}}-x)}\\
&+\frac{\frac{1}{n\alpha_{n}}\sum_{k=0}^{n-1}K_{h}(X_{t_{k}}-x)\left(\sigma\diamond(B_{t_{k+1}}-B_{t_{k}}^{H})\right)}{\frac{1}{n}\sum_{k=0}^{n-1}K_{h}(X_{t_{k}}-x)}\\
=&\frac{I+II+III}{S_{n,h}(x)}.\\
\end{aligned}
\end{equation}
We first state some lemmas in order to obtain our main consistency results.
\begin{lem}\label{lem1}
Under the conditions \ref{asp1}-\ref{asp5}, we have
\begin{equation*}
S_{n,h}(x)\rightarrow \mathbb{E}\big(K_{h}(\bar{X})\big),
\end{equation*}
almost surely, as $n\rightarrow\infty$.
\end{lem}
\noindent\textbf{Proof of Lemma \ref{lem1}}.
By properties of the Kernel function $K(\cdot)$, we have that there exists a constant $C_{K}<\infty$ and $p=1$ such that
\begin{equation*}
|K(x)|+|K^{\prime}(x)|\leq C_{K}\big|1+|x|\big|.
\end{equation*}
Then by Theorem \ref{th1}, we have
\begin{equation*}
\begin{aligned}
\frac{1}{n}\sum_{k=0}^{n-1}K_{h}(X_{t_{k}}-x)&=\frac{1}{n\alpha_{n}}\sum_{k=0}^{n-1}K_{h}(X_{t_{k}}-x)\alpha_{n}
=\frac{1}{t_{n}}\sum_{k=0}^{n-1}\int_{t_{k}}^{t_{k+1}}K_{h}(X_{t_{k}}-x)\mathrm{d}s\\
&=\frac{1}{t_{n}}\int_{0}^{t_{n}}\left\{\sum_{k=0}^{n-1}K_{h}(X_{t_{k}}-x)\mathbf{1}_{\left[t_{k}, t_{k+1}\right)}(s)\mathrm{d}s\right\}.
\end{aligned}
\end{equation*}
Therefore, the desired convergence result follows from Eq. (\ref{eq9}) directly.
\hfill$\square$

The following lemma shows that $I\rightarrow 0$ as $n\rightarrow \infty$.
\begin{lem}\label{lem2}
Under the conditions \ref{asp1}, \ref{asp2} and \ref{asp4}, we have 
\begin{equation*}
\sup_{t, s\in [t_{k},t_{k+1})}|X_{t}-X_{s}|\leq \left|\sigma\eta_{T}(t-s)^{\alpha}+b(X_{s})(t-s)\right|e^{L(t-s)},
\end{equation*}
almost surely, as $n\rightarrow\infty$.
\end{lem}
\noindent\textbf{Proof of Lemma \ref{lem2}}.
By using the Lipschitz condition of $b(\cdot)$ and Eq. (\ref{eq6}), we have
\begin{equation*}
\begin{aligned}
&\sup_{s\leq t\in [0,t_{n}]}|X_{t}-X_{s}|\\
&\leq\sup_{s\leq t\in [0,t_{n}]}\int_{s}^{t}b(X_{u})\mathrm{d}u+\sigma\eta_{T}|t-s|^{\alpha}\\
&=\sup_{s\leq t\in [0,t_{n}]}\int_{s}^{t}\left(b(X_{u})-b(X_{s})\right)\mathrm{d}u+b(X_{s})|t-s|+\sigma\eta_{T}|t-s|^{\alpha}\\
&\leq L\int_{s}^{t}\sup_{s\leq t\in [0,t_{n}]}|X_{t}-X_{s}|\mathrm{d}u+b(X_{s})|t-s|+\sigma\eta_{T}|t-s|^{\alpha}.
\end{aligned}
\end{equation*}
By Gronwall's inequality, we have
\begin{equation*}
\sup_{s\leq t\in [0,t_{n}]}|X_{t}-X_{s}|\leq \left|\sigma\eta_{T}(t-s)^{\alpha}+b(X_{s})(t-s)\right|e^{L(t-s)},
\end{equation*}
which completes the proof.
\hfill$\square$

Based on the Lemma \ref{lem2}, we could get the following lemma immediately.
\begin{lem}\label{lem3}
Under the conditions \ref{asp1}-\ref{asp5}, we have 
\begin{equation*}
\frac{1}{n\alpha_{n}}\sum_{k=0}^{n-1}K_{h}(X_{t_{k}}-x)\left(\int_{t_{k}}^{t_{k+1}}\left(b(X_{s})-b(X_{t_{k}})\mathrm{d}s\right)\right)\rightarrow 0 \quad \text{a.s.}.
\end{equation*}
\end{lem}
\noindent\textbf{Proof of Lemma \ref{lem3}}.
Based on the Lipschitz condition of $b(\cdot)$, we have
\begin{equation*}
\begin{aligned}
&\frac{1}{n\alpha_{n}}\sum_{k=0}^{n-1}K_{h}(X_{t_{k}}-x)\int_{t_{k}}^{t_{k+1}}\left(b(X_{s})-b(X_{t_{k}})\right)\mathrm{d}s\\
\leq&\frac{L}{n\alpha_{n}}\sum_{k=0}^{n-1}K_{h}(X_{t_{k}}-x)\int_{t_{k}}^{t_{k+1}}\left|X_{s}-X_{t_{k}}\right|\mathrm{d}s\\
\leq&\frac{L}{n}\sum_{k=0}^{n-1}K_{h}(X_{t_{k}}-x)\sup_{s\in [t_{k},t_{k+1})}\left|X_{s}-X_{t_{k}}\right|.
\end{aligned}
\end{equation*}
By using Lemma \ref{lem2}, $\alpha_{n}<1$ and $\alpha<H$, we have 
\begin{equation*}
\begin{aligned}
\sup_{s\in [t_{k},t_{k+1})}|X_{s}-X_{t_{k}}|
&\leq \left|\sigma\eta_{T}\alpha_{n}^{\alpha}+b(X_{s})\alpha_{n}\right|e^{L\alpha_{n}}\\
&\leq |\sigma\eta_{T}+b(X_{t_{k}})|\alpha_{n}^{\alpha}e^{L\alpha_{n}}.
\end{aligned}
\end{equation*}
By Theorem \ref{th1}, we have 
\begin{equation*}
\begin{aligned}
&\frac{L}{n}\sum_{k=0}^{n-1}K_{h}(X_{t_{k}}-x)\sup_{s\in [t_{k},t_{k+1})}\left|X_{s}-X_{t_{k}}\right|\\
\leq& LS_{n,h}(x)\sup_{s\in [t_{k},t_{k+1})}|X_{s}-X_{t_{k}}|\\
=&O(\alpha_{n}^{\alpha}),
\end{aligned}
\end{equation*}
which goes to $0$, as $n\rightarrow\infty$.
\hfill$\square$

Based on change of random variables and propertier of $K(\cdot)$, we give the following lemma.
\begin{lem}\label{lem4}
Under the conditions \ref{asp5} and \ref{asp6}, we have
\begin{equation*}
\frac{\frac{1}{n}\sum_{k=0}^{n-1}K_{h}(X_{t_{k}}-x)b(X_{t_{k}})}{\frac{1}{n}\sum_{k=0}^{n-1}K_{h}(X_{t_{k}}-x)}\rightarrow b(x) \quad \text{a.s.}.
\end{equation*}
\end{lem}
\noindent\textbf{Proof of Lemma \ref{lem4}}.
Denote $X_{t_{k}}=x+u_{t_{k}}h$, by using change of variables, we have
\begin{equation*}
\begin{aligned}
\frac{II}{S_{n,h}(x)}=&\frac{\frac{1}{n}\sum_{k=0}^{n-1}K_{h}(X_{t_{k}}-x)b(X_{t_{k}})}{\frac{1}{n}\sum_{k=0}^{n-1}K_{h}(X_{t_{k}}-x)}\\
=&\frac{\frac{1}{n}\sum_{k=0}^{n-1}K(u_{t_{k}})b(x+u_{t_{k}}h)}{\frac{1}{n}\sum_{k=0}^{n-1}K(u_{t_{k}})}.
\end{aligned}
\end{equation*} 
We do a further Taylor expansion of $\frac{II}{S_{n,h}(x)}$ around $h$ and obtain
\begin{equation*}
\begin{aligned}
\frac{II}{S_{n,h}(x)}=&\frac{\sum_{k=0}^{n-1}K(u_{t_{k}})\left(b(x)+hb^{\prime}(x+u_{t_{k}}h)+o(h^{2}))\right)}{\sum_{k=0}^{n-1}K(u_{t_{k}})}\\
=&b(x)+h\frac{\sum_{k=0}^{n-1}K(u_{t_{k}})b^{\prime}(x+u_{t_{k}}h)}{\sum_{k=0}^{n-1}K(u_{t_{k}})}+o(h^{2})\\
\rightarrow& b(x)  \quad a.s.,
\end{aligned}
\end{equation*}
which completes the proof.
\hfill$\square$

The following lemma shows that $III\rightarrow 0$, as $n\rightarrow\infty$.
\begin{lem}\label{lem5}
	Under conditions \ref{asp4}, \ref{asp5} and \ref{asp6}, we have
	\begin{equation*}
		\frac{1}{n\alpha_{n}}\sum_{k=0}^{n-1}K_{h}(X_{t_{k}}-x)\left(\sigma \diamond(B_{t_{k+1}}-B_{t_{k}}^{H})\right)\rightarrow 0,
	\end{equation*}
in probability, as $n\rightarrow\infty$.
\end{lem}
\noindent\textbf{Proof of Lemma \ref{lem5}}. 
Note that the It$\hat{o}$-Skorokhod integral coulded be defined as the Riemann sum using Wick product Eq. (\ref{eq8}), then
\begin{equation*}
	\lim_{n\rightarrow\infty}III=\frac{\sigma}{t_{n}}\int_{0}^{t_{n}}K_{h}(X_{t_{k}}-x)\mathbf{1}_{\left[t_{k}, t_{k+1}\right)}(t)\mathrm{d}B_{t}^{H}.
\end{equation*}
From the properties of It$\hat{o}$ integral, we have 
\begin{equation*}
\mathbb{E}\big(\int_{0}^{t_{n}}K_{h}(X_{t_{k}}-x)\mathbf{1}_{\left[t_{k}, t_{k+1}\right)}(t)\mathrm{d}B^{H}_{t}\big)=0.
\end{equation*}
As a consequence of Meyer's inequality \citep{Hu(2005)}, we have
\begin{equation*}
	\begin{aligned}
		&\mathbb{E}\bigg(\int_{0}^{t_{n}}K_{h}(X_{t_{k}}-x)\mathbf{1}_{\left[t_{k}, t_{k+1}\right)}(t)\mathrm{d}B^{H}_{t}\bigg)^{2}\\
		\leq&C_{p}\big(\mathbb{E}\|K_{h}\|^{2}_{\mathcal{H}}+\mathbb{E}\|DK_{h}\|^{2}_{\mathcal{H}\otimes\mathcal{H}}\big),
	\end{aligned}
\end{equation*}
where $C_{p}$ is a constant.
By Eq. (\ref{eq7}), we have
\begin{equation*}
\begin{aligned}
&\mathbb{E}\bigg(\int_{0}^{t_{n}}K_{h}(X_{t_{k}}-x)\mathbf{1}_{\left[t_{k}, t_{k+1}\right)}(t)\mathrm{d}B^{H}_{t}\bigg)^{2}\\
\leq&C_{p,H}\bigg[\bigg(\int_{0}^{t_{n}}\mathbb{E}\big|K_{h}(X_{t_{k}}-x)\mathbf{1}_{[t_{k}, t_{k+1})}(t)\big|^{\frac{1}{H}}\mathrm{d}t\bigg)^{2H}\\
&+\mathbb{E}\bigg(\int_{0}^{T}\int_{0}^{T}\big|K^{\prime}_{h}(X_{t_{k}}-x)\mathbf{1}_{[t_{k}, t_{k+1})}(t)\mathrm{D}_{s}X_{t_{k}}\big|^{\frac{1}{H}}\mathrm{d}s\mathrm{d}t\bigg)^{2H}\bigg].
\end{aligned}
\end{equation*}
Under conditions \ref{asp5}, we have
\begin{equation*}
\begin{aligned}
\bigg(\int_{0}^{t_{n}}\mathbb{E}\big|K_{h}(X_{t_{k}}-x)\mathbf{1}_{[t_{k}, t_{k+1})(t)}\big|^{\frac{1}{H}}\mathrm{d}t\bigg)^{2H}\leq \|K_{h}\|_{\infty}t_{n}^{2H}=t^{2H}h^{-2}\|K\|_{\infty}.
\end{aligned}
\end{equation*}
Under conditions \ref{asp6} and Proposition \ref{prop1}, we have
\begin{equation*}
\begin{aligned}
&\bigg(\int_{0}^{t_{n}}\int_{0}^{t}\big|K^{\prime}_{h}(X_{t_{k}}-x)\mathbf{1}_{[t_{k}, t_{k+1})}(t)\mathrm{D}_{s}X_{t_{k}}\big|^{\frac{1}{H}}\mathrm{d}s\mathrm{d}t\bigg)^{2H}\\
\leq&\bigg(\int_{0}^{t_{n}}\int_{0}^{t}\big|K^{\prime}_{h}(X_{t_{k}}-x)\mathbf{1}_{[t_{k}, t_{k+1})}(t)\sigma e^{-M(t_{k}-s)}\big|^{\frac{1}{H}}\mathrm{d}s\mathrm{d}t\bigg)^{2H}
\\
\leq&\frac{1}{M^{2H}}\|K_{h}^{\prime}\|_{\infty}^{2}\sigma^{2}\bigg(\int_{0}^{t_{n}}(1-e^{-\frac{Mt}{H}})\mathrm{d}t\bigg)^{2H}\\
=&O(h^{-4}t_{n}^{2H}).
\end{aligned}
\end{equation*}
Combining the above two equations, we have
\begin{equation*}
\mathbb{E}\bigg(\int_{0}^{t_{n}}K\big(\frac{X_{t_{k}}-x}{h}\big)\mathbf{1}_{\left[t_{k}, t_{k+1}\right)}(t)\mathrm{d}B^{H}_{t}\bigg)^{2}\leq O(h^{-4}t_{n}^{2H})
\end{equation*}
By Chebyshev's Inequality, we have that for $\gamma\in(H-1,0)$,
\begin{equation*}
\begin{aligned}
&\mathbb{P}\bigg(\frac{1}{t_{n}}\sum_{k=0}^{n-1}K_{h}(X_{t_{k}}-x)\left(\sigma\diamond(B_{t_{k+1}}-B_{t_{k}}^{H})\right)\geq t_{n}^{\gamma}h^{-2}\bigg)\leq O(t_{n}^{2H-2-2\gamma}),
\end{aligned}
\end{equation*}
which goes to $0$ as $n\rightarrow\infty$.
\hfill$\square$

The following theorem establishes the consistency of this estimator.
\begin{thm}\label{th2}
Under the conditions \ref{asp1}-\ref{asp6}, we have that the estimator $\hat{b}_{n,h}(x)$ defined in Eq. (\ref{eq5}) convergences to $b(x)$ in probability, as $n\rightarrow\infty$.
\end{thm}
\noindent\textbf{Proof of Theorem \ref{th2}}.
For $H\in(\frac{1}{2},1)$, we have that the estimator defined in Eq. (\ref{eq5}) is equivalent to the estimator defined in Eq. (\ref{eq10}).
By Lemma \ref{lem1} and Lemma \ref{lem3}, we have 
\begin{equation*}
\frac{I}{S_{n,h}(x)}\rightarrow 0 \quad \text{a.s.}.
\end{equation*}
By Lemma \ref{lem1} and Lemma \ref{lem4}, we have 
\begin{equation*}
\frac{II}{S_{n,h}(x)}\rightarrow b(x) \quad \text{a.s.}.
\end{equation*}
By Lemma \ref{lem1} and \ref{lem5}, we have
\begin{equation*}
\frac{III}{S_{n,h}(x)}\rightarrow 0 \quad \text{in probability}.
\end{equation*}
Combining the above three equations, we obtain the desired results.
\hfill$\square$

\end{document}